\documentclass[12pt]{article}

\usepackage{amsmath,amssymb}
\usepackage[applemac]{inputenc}
\usepackage{amsmath,amssymb,fullpage}
\usepackage{graphicx}
\usepackage{showlabels}
\usepackage{color}

\newcommand{\dproof}{\noindent {Proof.} \quad}
\newcommand{\fproof}{\hfill $\square$ \bigskip}

\newtheorem{definition}{Definition}[section]
\newtheorem{example}{Example}[section]
\newtheorem{theorem}[definition]{Theorem}
\newtheorem{problem}[definition]{Problem}
\newtheorem{remark}[definition]{ \it Remark}

\numberwithin{equation}{section}

\def\1B{\text{1\!\!I}}

\begin{document}
\date{30 August 2019 }
\title{A new approach to optimal stopping for Hunt processes}
\author{
Achref Bachouch$^{1,2}$, Olfa Draouil$^{1,2}$ and Bernt \O ksendal$^{1,2}$}

\footnotetext[1]{Department of Mathematics, University of Oslo, P.O. Box 1053 Blindern, N--0316 Oslo, Norway.\\
Emails: {\tt achrefb@math.uio.no, \tt olfad@math.uio.no, \tt oksendal@math.uio.no}}

\footnotetext[2]{This research was carried out with support of the Norwegian Research Council, within the research project Challenges in Stochastic Control, Information and Applications (STOCONINF), project number 250768/F20.}

\maketitle


\begin{abstract}
In this paper we present a new verification theorem for optimal stopping problems for Hunt processes. The approach is based on the Fukushima-Dynkin formula \cite{F} and its advantage is that it allows us to verify that a given function is the value function without using the viscosity solution argument. Our verification theorem works in any dimension.
We illustrate our results with some examples of optimal stopping of reflected diffusions and absorbed diffusions.
\end{abstract}

\section{Introduction}
Usually when solving optimal stopping problems it is assumed that the value function is twice continuously differentiable  ($C^2$), and in order to find it we use the high contact principle.
But in reality this is a strong condition and sometimes the high contact principle is not valid. Then one must use the viscosity solution approach to verify that a given function is indeed the value function.
This is done in for example Dai \& Menoukeu Pamen \cite{DM}, where  the authors use a viscosity solution approach to study optimal stopping of some processes with reflection.
More precisely, they prove that the value function is the unique viscosity solution of the HJB equation associated with the optimal stopping problem of reflected Feller processes.\\

 In our paper we treat a more general case by considering the Hunt processes, i.e. strong Markov and quasi left continuous processes with respect a filtration $\{\mathcal{F}_t\}_{t\geq 0}$ the natural filtration of the Hunt process. We come back later in the next section with more details about Hunt processes.
Using the Fukushima-Dynkin formula \cite{F}, we obtain a variational inequality verification theorem for optimal stopping of Hunt processes.
Note that in this theorem we do not need the use of high contact principle, nor do we need the viscosity approach. Moreover, our theorem works for multidimensional case, i.e. our Hunt process $X$ can take value in $\mathbb{R}^d$ for any $d\geq 1$.Therefore our work can be regarded as a generalisation of the result in the Mordecki \& Salminen paper \cite{Salminen}, which is based on methods applicable only in the 1-dimensional case.\\

This paper is organised as follow:\\
In Section 2 we define the Hunt process, we present the Fukushima-Dynkin formula \cite{F} then we make a connection between Fukushima-Dynkin formula  and the Dynkin formula in \O ksendal-Sulem formula \cite{SO}.\\

In Section 3 we introduce the optimal stopping problem for Hunt processes. We prove a verification theorem for this problem.\\

Finally in Section 4 we illustrate our method by applying it to some optimal stopping problems for reflected or absorbed Hunt processes.

\section{Our main result: \\
Variational inequalities for optimal stopping}

In the following we let 
$\{X_t\}_{0\leq t\leq T}$  be a given Hunt process with values in a closed subset $K$ of $\mathbb{R}^d$. Our method and results are valid in any dimension, but for simplicity of notation we will in the following assume that $d=1.$ We let $\mathcal{T}$ be the set of stopping times $\tau \leq T$ with respect to the filtration $\mathbb{F}=\{\mathcal{F}_t\}_{t\geq 0}$ generated by $X$. We refer to the book by Fukushima \cite{F} for more information about Hunt processes and their associated Dirichlet forms and calculus.

We consider the following optimal stopping problem: 
\begin{problem}Given functions $f$ and $g$ and a constant $\alpha > 0$, find $\Phi_{\alpha}(x)$ and $\tau^* \in \mathcal{T}$ such that
\begin{equation}
\Phi_{\alpha}(x)=\sup_{\tau\in \mathcal{T}}J_{\alpha}^{\tau}(x)=J_{\alpha}^{\tau^*}(x)
\end{equation}
where
\begin{equation}
J_{\alpha}^{\tau}(x)=E_x[\int_0^{\tau}e^{-\alpha t}f(X_t)dt+e^{-\alpha \tau}g(X_{\tau})]; \quad x \in K.
\end{equation}

\end{problem}

To study this problem we introduce the \emph{resolvent} of $X$, defined by
\begin{equation}
R_{\alpha}\varphi(x) = E_{x}[\int_0^{\infty} e^{-\alpha t} \varphi(X(t)) dt]
\end{equation}
for all functions $\varphi$ such that the integral converges.

Recall that the \emph{Fukushima-Dynkin formula} (see (4.2.6) p. 97 in \cite{F}), can be written
\begin{equation}
E_x[e^{-\alpha \tau} R_{\alpha} \psi (X_{\tau)}]=R_{\alpha}\psi(x)-E_x[\int_0^{\tau} e^{-\alpha t}\psi(X_t)dt];\quad \text{ for all stopping times } \tau.
\end{equation}

Using this, we obtain the following result, which is our main result. It can be regarded as a weak analogue of the variational inequality (4) of \cite{P}. Note that a viscosity solution interpretation is not needed here:

\begin{theorem}(Variational inequalities for optimal stopping)

Suppose there exists a measurable function $\psi:K\mapsto \mathbb{R}$ such that thew following hold:\\
(i) $R_{\alpha}\psi(x) \geq g(x)$ for all $x\in K,$\\
(ii) $\psi(x) \geq f(x)$ for all $x\in K,$\\
(iii) Define $D=\{ x\in K, R_{\alpha}\psi(x)>g(x)\}$,\\
and put\\
$\tau_{D}=\inf\{ t>0, X_t\notin D\}.$\\
Assume that \\
(iv) $\psi(x)=f(x)$ on $D$,\\
(v) the family $\{R_{\alpha}\psi(X_{\tau}), \tau \leq \tau_D \}$ is uniformly integrable with respect to $P_x$ for all $x\in K$.\\
Then
\begin{equation}
R_{\alpha}\psi(x)=\Phi_{\alpha}(x)=\sup_{\tau }J_{\alpha}(x), \forall x\in K
\end{equation}
and
\begin{equation}
\tau^*=\tau_D
\end{equation}
is an optimal stopping time.
\end{theorem}

\dproof
First we remark that for $\tau=0$ we have
\begin{equation}\label{tau0}
J_{\alpha}^0(x)=g(x)\leq \Phi_{\alpha}(x).
\end{equation}
Using (i),(ii) and the Fukushima-Dynkin formula we get that
\begin{align}
R_{\alpha}\psi(x)&=E_x[\int_0^{\tau} e^{-\alpha t}\psi(X_t)dt]+E_x[e^{-\alpha \tau} R_{\alpha} \psi (X_{\tau)}]\nonumber\\
&\geq E_x[\int_0^{\tau}e^{-\alpha t}f(X_t)dt]+E_x[e^{-\alpha \tau}g(X_{\tau})]=J_{\alpha}^{\tau}(x), \forall x\in X.
\end{align}
Since this inequality holds for arbitrary $\tau\leq T$ then we get
\begin{equation}\label{ineq1}
R_{\alpha}\psi(x)\geq \Phi_{\alpha}(x),\forall x\in K.
\end{equation}
To prove the reverse inequality we consider two cases
\begin{itemize}
\item Suppose that $x\notin D$ then
\begin{equation}\label{xnotinD}
R_{\alpha}\psi(x)=g(x)\leq \Phi_{\alpha}(x).
\end{equation}
Note that the last inequality of equation \eqref{xnotinD} comes from \eqref{tau0}.
Combining \eqref{ineq1} and \eqref{xnotinD} we conclude that
\begin{equation}
J_{\alpha}^0(x)=R_{\alpha}\psi(x)=\Phi_{\alpha}(x),\forall x\notin D \text{ and } \tau^*=\tau^*(x,\omega)=0.
\end{equation}

\item Suppose $x\in D$. Let $\{D_k\}_{k\geq 1}$ be a sequence of increasing open sets of $D_k$ such that $\bar{D_k}\subset D$, $\bar{D_k}$ is compact and $D=\cup_{k=1}^{\infty}D_k$. \\
Let $\tau_k=\inf\{ t>0, X_t\notin D_k\}$.\\
Choose $x\in D_k$. Then by the Fukushima-Dynkin formula and (iv) we have
\begin{align}
R_{\alpha}\psi(x)&=E_x[\int_0^{\tau_k} e^{-\alpha t}\psi(X_t)dt+e^{-\alpha \tau_k} R_{\alpha} \psi (X_{\tau_k)}]\nonumber\\
&=E_x[\int_0^{\tau_k} e^{-\alpha t}f(X_t)dt+e^{-\alpha \tau_k} R_{\alpha} \psi (X_{\tau_k)}].
\end{align}
Using the uniform integrability of $R_{\alpha}\psi(X_{\tau_k})$, $\tau_k<\tau_D$, quasi left-continuity and the fact that $R_{\alpha}\psi(X_{\tau_D})=g(X_{\tau_D})$ we get
\begin{align}\label{xinD}
R_{\alpha}\psi(x)&=\lim_{k\rightarrow +\infty} E_x[\int_0^{\tau_k} e^{-\alpha t}f(X_t)dt+e^{-\alpha \tau_k} R_{\alpha} \psi (X_{\tau_k)}]\\
&=E_x[\int_0^{\tau_D} e^{-\alpha t}f(X_t)dt+e^{-\alpha \tau_D} R_{\alpha} \psi (X_{\tau_D)}]\\
&=E_x[\int_0^{\tau_D} e^{-\alpha t}f(X_t)dt+e^{-\alpha \tau_D} g(X_{\tau_D)}]=J_{\alpha}^{\tau_D}(x)\leq \Phi_{\alpha}(x).
\end{align}
Combining \eqref{ineq1} and \eqref{xinD} we get
\begin{equation}
\Phi_{\alpha}(x)\leq R_{\alpha}\psi(x)=J_{\alpha}^{\tau_D}(x)\leq \Phi_{\alpha}(x).
\end{equation}
Then we conclude that
\begin{equation}
R_{\alpha}\psi(x)=\Phi_{\alpha}(x) \text{ and } \tau^*=\tau^*(x,\omega)=\tau_D, \forall x\in D.
\end{equation}
\end{itemize}
\fproof

\section{Examples}
To  illustrate our main result, we study some examples:

\begin{example}
Consider the following optimal stopping problem:\\

Find $\Phi_{\alpha}(x)$ and $\tau^*$ such that
\begin{equation}
\Phi_{\alpha}(x)= \sup_{\tau} E[\int_0^{\tau} e^{-\alpha t}B(t)dt]=E[\int_0^{\tau*}e^{-\alpha t}B(t)dt].
\end{equation}

We want to solve this problem in two ways:\\
(i)  By using the classical variational inequality theorem approach in the SDE book \\
(ii) By using Theorem 0.1 above.

\end{example}
(i)\\
Gessing $D=\{(t,x), x> x_0\}, x_0<0$.\\
The function $\phi_{\alpha}$ should verify the following PDE
\begin{equation}
\begin{cases}
\frac{\partial \phi_{\alpha}}{\partial t}(t,x)+\frac{1}{2}\frac{\partial^2 \phi_{\alpha}}{\partial x^2}(t,x)+xe^{-\alpha t}=0 \text{ on } D\\
\phi_{\alpha}(t,x)=0,\quad x\notin D.
\end{cases}
\end{equation}
Put $\phi_{\alpha}(t,x)=\phi_0(x) e^{-\alpha t}$
then $\phi_0$ verifies the following second order differential equation
\begin{equation}\label{secondordereq}
\begin{cases}
\frac{1}{2}\frac{\partial^2 \phi_{0}}{\partial x^2}(x)-\alpha \phi_{0}(x) +x=0 \text{ on } D\\
\phi_{0}(t,x)=0,\quad x\notin D.
\end{cases}
\end{equation}
The general solution of the equation \eqref{secondordereq} is given by
\begin{equation}
\phi_0(x)=C_1e^{\sqrt{2\alpha}x}+C_2e^{-\sqrt{2\alpha}x}
+\frac{1}{\alpha}x
\end{equation}
where $C_1$ and $C_2$ are constants.\\
Since we have
\begin{equation}
E_x[\int_0^{\tau} e^{-\alpha t} B(t)dt] = \frac{1}{\alpha}(1-e^{-\alpha \tau})x + E_0[\int_0^{\tau} e^{-\alpha t} B(t)dt]
\end{equation}

then only the first term depends on $x$ and it grows at most linearly as $x$ goes to $+\infty$.
Hence $C_1=0$.

Using the continuity of $\phi_{\alpha}$ at $x=x_0$ we have $\phi_0(x_0)=0$ then
\begin{equation}
C_2e^{-\sqrt{2\alpha}x_0}+\frac{1}{\alpha}x_0=0.
\end{equation}

 Hence $C_2=-\frac{1}{\alpha}x_0 e^{\sqrt{2\alpha}x_0}$.

Then
\begin{equation}
\phi_{\alpha}(t,x)=e^{-\alpha t}(-\frac{1}{\alpha}x_0 e^{\sqrt{2\alpha}x_0}e^{-\sqrt{2\alpha}x}
+\frac{1}{\alpha}x).
\end{equation}

Using now the high contact equation i.e,  $\phi_{\alpha}$ is $C^1$ at $x=x_0$ we get the following equation
\begin{equation}
-\frac{1}{\alpha}e^{\sqrt{2\alpha}x_0}e^{-\sqrt{2\alpha}x_0}(-\sqrt{2\alpha})+\frac{1}{\alpha}=0
\end{equation}
Then we deduce that
\begin{equation}
x_0 =-\frac{1}{\sqrt{2\alpha}}.
\end{equation}

(ii)We now solve the problem using our approach.
By condition (iii) of Theorem 2.2 we have $\psi(x)=f(x)$ on $D=]x_0,+\infty[$. i:e $\psi(x)=x$.
\begin{align}\label{ex1R}
R_{\alpha}\psi(x)&= \int_{x_0}^{+\infty}\psi(y)R_{\alpha}(x,dy)\\
&=\int_{x_0}^{+\infty}\psi(y)E_x[\int_{0}^{+\infty}e^{-\alpha t}1_{dy}(B_t)dt]\\
&=\int_{0}^{+\infty}e^{-\alpha t}\int_{x_0}^{+\infty}yE_x[1_{dy}(B_t)]dt\\
&=\int_{0}^{+\infty}e^{-\alpha t}\int_{x_0}^{+\infty}yP_x(B_t\in dy)dt\label{eq2.31}\\
&=\int_{0}^{+\infty}e^{-\alpha t}\int_{x_0}^{+\infty}y\frac{e^{-\frac{(y-x)^2}{2t}}}{\sqrt{2\pi t}}dy dt\\
&=\int_{0}^{+\infty}e^{-\alpha t}\int_{x_0-x}^{+\infty}(z+x)\frac{e^{-\frac{z^2}{2t}}}{\sqrt{2\pi t}}dz dt\\
&=\int_{0}^{+\infty}e^{-\alpha t}\int_{x_0-x}^{+\infty}z\frac{e^{-\frac{z^2}{2t}}}{\sqrt{2\pi t}}dz dt+x\int_{0}^{+\infty}e^{-\alpha t}\int_{x_0-x}^{+\infty}\frac{e^{-\frac{z^2}{2t}}}{\sqrt{2\pi t}}dz dt                                                                                                                                                                                                                                                                                                                                                                                                                                                                                                                                                                                                                                                                                                                                                                                                                                                                                                                                                                                                                                                                                                                                                                                                                                                                                                                                                                                                                                                                                                                                                                                                                                                                                                                                                                                                                                                                                                                                                                                                                                                                                                                                                                                                                                                                                                                                                                                                                                                                                                                                                                                                                                                                                                                                                                                                                                                                                                                                                                                                                                                                                                                                                                                                                                                                                                                                                                                                                                                                                                                                                                                                                                                                                                                                                                                                                                                              \\
&=\int_{0}^{+\infty}e^{-\alpha t}t\frac{e^{-\frac{(x_0-x)^2}{2t}}}{\sqrt{2\pi t}}dt+
x\int_{0}^{+\infty}e^{-\alpha t}\int_{x_0-x}^{+\infty}\frac{e^{-\frac{z^2}{2t}}}{\sqrt{2\pi t}}dz dt\label{eq3.17}
\end{align}

We distinguish two cases:

\begin{itemize}

\item 1) If $x>x_0$ then $x_0-x<0$.
In this case we have from equation \eqref{eq3.17} that 
\begin{align}
R_{\alpha}\psi(x)&=\int_{0}^{+\infty}e^{-\alpha t}t\frac{e^{-\frac{(x_0-x)^2}{2t}}}{\sqrt{2\pi t}}dt+
x\int_{0}^{+\infty}e^{-\alpha t}\int_{x_0-x}^{+\infty}\frac{e^{-\frac{z^2}{2t}}}{\sqrt{2\pi t}}dz dt\\
&=\int_{0}^{+\infty}e^{-\alpha t}t\frac{e^{-\frac{(x_0-x)^2}{2t}}}{\sqrt{2\pi t}}dt+
x\int_{0}^{+\infty}e^{-\alpha t}\int_{x_0-x}^{0}\frac{e^{-\frac{z^2}{2t}}}{\sqrt{2\pi t}}dz dt+
x\int_{0}^{+\infty}e^{-\alpha t}\int_{0}^{+\infty}\frac{e^{-\frac{z^2}{2t}}}{\sqrt{2\pi t}}dz dt\\
&=\int_{0}^{+\infty}e^{-\alpha t}t\frac{e^{-\frac{(x_0-x)^2}{2t}}}{\sqrt{2\pi t}}dt+x\int_{0}^{+\infty}e^{-\alpha t}\int_0^{x-x_0}\frac{e^{-\frac{z^2}{2t}}}{\sqrt{2\pi t}}dz dt+
\frac{x}{2}\int_{0}^{+\infty}e^{-\alpha t}dt\\
&=\int_{0}^{+\infty}e^{-\alpha t}t\frac{e^{-\frac{(x_0-x)^2}{2t}}}{\sqrt{2\pi t}}dt+x\int_{0}^{+\infty}e^{-\alpha t}\int_0^{\frac{x-x_0}{\sqrt{2t}}}\frac{e^{-z^2}}{\sqrt{\pi }}dz dt+
\frac{x}{2}\int_{0}^{+\infty}e^{-\alpha t}dt\\
&=\int_{0}^{+\infty}e^{-\alpha t}t\frac{e^{-\frac{(x_0-x)^2}{2t}}}{\sqrt{2\pi t}}dt+\frac{x}{2}\int_{0}^{+\infty}e^{-\alpha t}(\int_0^{\frac{x-x_0}{\sqrt{2t}}}\frac{e^{-z^2}}{\sqrt{\pi }}dz+1)dt
\\
&=\int_{0}^{+\infty}e^{-\alpha t}t\frac{e^{-\frac{(x_0-x)^2}{2t}}}{\sqrt{2\pi t}}dt+\frac{x}{2}\int_{0}^{+\infty}e^{-\alpha t}(\int_0^{\frac{x-x_0}{\sqrt{2t}}}\frac{2e^{-z^2}}{\sqrt{\pi }}dz-1)dt+x\int_{0}^{+\infty}e^{-\alpha t}dt\\
&=\int_{0}^{+\infty}e^{-\alpha t}t\frac{e^{-\frac{(x_0-x)^2}{2t}}}{\sqrt{2\pi t}}dt+\frac{x}{\alpha}+\frac{x}{2}\int_{0}^{+\infty}e^{-\alpha t}(\int_0^{\frac{x-x_0}{\sqrt{2t}}}\frac{2e^{-z^2}}{\sqrt{\pi }}dz-1)dt\label{eq3.29}\\
&=\int_{0}^{+\infty}e^{-\alpha t}t\frac{e^{-\frac{(x_0-x)^2}{2t}}}{\sqrt{2\pi t}}dt+\frac{x}{\alpha}-\frac{x}{2}\int_{0}^{+\infty}e^{-\alpha t}(1-H(\frac{x-x_0}{\sqrt{2t}}))dt
\end{align}
where 
\begin{equation}
H(x)=\int_0^{x}\frac{2e^{-z^2}}{\sqrt{\pi }}dz.
\end{equation}
By \cite{GR} we have the following result
\begin{equation}
\int_0^{+\infty}(1-H(\frac{q}{2\sqrt{t}}))e^{-\alpha t} dt=\frac{1}{\alpha}e^{-q\sqrt{\alpha}}, \quad Re\text{ } \alpha>0, |arg \text{ }q|<\frac{\pi}{4}.
\end{equation}
Then, for $q=\sqrt{2}(x-x_0)$ we get 
\begin{equation}\label{eq3.33}
\int_{0}^{+\infty}e^{-\alpha t}(1-H(\frac{x-x_0}{\sqrt{2t}}))dt=\frac{1}{\alpha}e^{-\sqrt{2\alpha}(x-x_0)}
\end{equation}
Therefore
\begin{equation}\label{eq3.34}
R_{\alpha}\psi(x) =\int_{0}^{+\infty}e^{-\alpha t}t\frac{e^{-\frac{(x_0-x)^2}{2t}}}{\sqrt{2\pi t}}dt+\frac{x}{\alpha}-\frac{x}{2\alpha}e^{-\sqrt{2\alpha}(x-x_0)}.
\end{equation}
From \eqref{eq3.33} we have 
\begin{equation}\label{eq3.35}
\int_{0}^{+\infty}e^{-\alpha t}(1-\int_0^{\frac{x-x_0}{\sqrt{2t}}}\frac{2e^{-z^2}}{\sqrt{\pi }}dz)dt=\frac{1}{\alpha}e^{-\sqrt{2\alpha}(x-x_0)}.
\end{equation}
Derive equation \eqref{eq3.35} with respect to $x$ we find:
\begin{equation}\label{eq3.36}
2\int_{0}^{+\infty}e^{-\alpha t}\frac{e^{-\frac{(x-x_0)^2}{2t}}}{\sqrt{2\pi t}}dt=\sqrt{\frac{2}{\alpha}}e^{-\sqrt{2\alpha}(x-x_0)}.
\end{equation}
Derive equation \eqref{eq3.36} with respect to $\alpha$ we get 
\begin{equation}\label{eq3.37}
\int_{0}^{+\infty}te^{-\alpha t}\frac{e^{-\frac{(x-x_0)^2}{2t}}}{\sqrt{2\pi t}}dt=e^{-\sqrt{2\alpha}(x-x_0)}(\frac{1}{\alpha\sqrt{2\alpha}}+\frac{x}{2\alpha}).
\end{equation}
Replacing equation\eqref{eq3.37} in \eqref{eq3.34} we get 
\begin{equation}\label{eq3.38b}
R_{\alpha}\psi(x)=\frac{1}{\alpha\sqrt{2\alpha}}e^{-\sqrt{2\alpha}(x-x_0)}+\frac{x}{\alpha}>0, \quad x>x_0.
\end{equation}
\item 2) If $x\leq x_0$ then $x_0-x\geq 0$, then equation \eqref{eq3.17} becomes
\begin{align}
&R_{\alpha}\psi(x)=\int_{0}^{+\infty}e^{-\alpha t}t\frac{e^{-\frac{(x_0-x)^2}{2t}}}{\sqrt{2\pi t}}dt+
x\int_{0}^{+\infty}e^{-\alpha t}\int_{x_0-x}^{+\infty}\frac{e^{-\frac{z^2}{2t}}}{\sqrt{2\pi t}}dz dt\\
&=\int_{0}^{+\infty}e^{-\alpha t}t\frac{e^{-\frac{(x_0-x)^2}{2t}}}{\sqrt{2\pi t}}dt+
x\int_{0}^{+\infty}e^{-\alpha t}(\int_{0}^{+\infty}\frac{e^{-\frac{z^2}{2t}}}{\sqrt{2\pi t}}dz -\int_0^{x_0-x}\frac{e^{-\frac{z^2}{2t}}}{\sqrt{2\pi t}}dz)dt\\
&=\int_{0}^{+\infty}e^{-\alpha t}t\frac{e^{-\frac{(x_0-x)^2}{2t}}}{\sqrt{2\pi t}}dt+
\frac{x}{2}\int_{0}^{+\infty}e^{-\alpha t}(1-\int_0^{\frac{x_0-x}{\sqrt{2t}}}\frac{2e^{-z^2}}{\sqrt{\pi}}dz)dt.
\end{align}
Or we have

\begin{equation}\label{eq2.50}
\int_{0}^{+\infty}e^{-\alpha t}(1-\int_0^{\frac{x_0-x}{\sqrt{2t}}}\frac{2 e^{-z^2}}{\sqrt{\pi}}dz)dt=\frac{1}{\alpha}e^{-\sqrt{2\alpha}(x_0-x)}
\end{equation}
Then
\begin{equation}\label{eq2.51}
R_{\alpha}\psi(x)=\int_{0}^{+\infty}e^{-\alpha t}t\frac{e^{-\frac{(x_0-x)^2}{2t}}}{\sqrt{2\pi t}}dt+\frac{x}{2\alpha}e^{-\sqrt{2\alpha}(x_0-x)}.
\end{equation}
Let us derive \eqref{eq2.50} with respect to $x$, then we get
\begin{equation}\label{eq2.52}
2\int_{0}^{+\infty}e^{-\alpha t}\frac{e^{-\frac{(x_0-x)^2}{2t}}}{\sqrt{2\pi t}}dt=\frac{\sqrt{2}}{\sqrt{\alpha}}e^{-\sqrt{2\alpha}(x_0-x)}.
\end{equation}
Let us now derive \eqref{eq2.52} with respect to $\alpha$ then we get:
\begin{equation}
-2\int_{0}^{+\infty}te^{-\alpha t}\frac{e^{-\frac{(x_0-x)^2}{2t}}}{\sqrt{2\pi t}}dt=
e^{-\sqrt{2\alpha}(x_0-x)}\sqrt{2}(-\frac{1}{2\alpha\sqrt{\alpha}}-\frac{\sqrt{2}(x_0-x)}{2\alpha}).
\end{equation}
Then
\begin{equation}\label{eq3.46b}
\int_{0}^{+\infty}te^{-\alpha t}\frac{e^{-\frac{(x_0-x)^2}{2t}}}{\sqrt{2\pi t}}dt=
e^{-\sqrt{2\alpha}(x_0-x)}(\frac{1}{2\alpha\sqrt{2\alpha}}+\frac{x_0-x}{2\alpha}).
\end{equation}
Substituting \eqref{eq3.46b} into \eqref{eq2.51} we get
\begin{equation}\label{eq3.47b}
R_{\alpha}\psi(x)=e^{-\sqrt{2\alpha}(x_0-x)}(\frac{1}{2\alpha\sqrt{2\alpha}}+\frac{x_0}{2\alpha}).
\end{equation} 
By continuity of $R_{\alpha}\psi$ at $x=x_0$ and combining \eqref{eq3.47b} with \eqref{eq3.38b} at $x=x_0$ we get 
\begin{equation}
x_0=-\frac{1}{\sqrt{2\alpha}}.
\end{equation}
Replacing the value of $x_0$ in \eqref{eq3.47b} we get that 
\begin{equation}
R_{\alpha}\psi(x)=0, \quad x\leq x_0.
\end{equation}

\end{itemize}
With $x_0=-\frac{1}{\sqrt{2\alpha}}$ it is easy to check that  $R_{\alpha}\psi(x)>0$ for  $x\in D=\{x>x_0\}$.\\

Next we consider the following example:

\begin{example}
Find $\Phi_{\alpha}(x)$ and $\tau^{*}$ such that
\begin{equation}
\Phi_{\alpha}(x)= \sup_{\tau} E[e^{-\alpha \tau} B(\tau)] = E[e^{-\alpha \tau^{*}} B(\tau^*)].
\end{equation}
Again we want to solve this problem in two ways:\\
(i)  By using the classical variational inequality theorem approach in the SDE book \\
(ii) By using Theorem 0.1 above.
\end{example}
By waiting long enough we can always get a payoff which is positive, because the exponential goes to 0 faster than the Brownian motion grows, eventually. Therefore it does not make sense to stop while $B(t) < 0$. So the continuation region should be of the form:
\begin{equation}
D=\{x<x_0\}, x_0>0.
\end{equation}

The function $\phi_{\alpha}$ should verify the following PDE
\begin{equation}
\begin{cases}
\frac{\partial \phi_{\alpha}}{\partial t}(t,x)+\frac{1}{2}\frac{\partial^2 \phi_{\alpha}}{\partial x^2}(t,x)=0 \text{ on } D\\
\phi_{\alpha}(t,x)=e^{-\alpha t} x,\quad x\notin D.
\end{cases}
\end{equation}
Put $\phi_{\alpha}(t,x)=\phi_0(x) e^{-\alpha t}$
then $\phi_0$ verifies the following second order differential equation
\begin{equation}\label{secondordereq1}
\begin{cases}
\frac{1}{2}\frac{\partial^2 \phi_{0}}{\partial x^2}(x)-\alpha \phi_{0}(x)=0 \text{ on } D\\
\phi_{0}(x)=x,\quad x\notin D.
\end{cases}
\end{equation}
The general solution of the equation \eqref{secondordereq1} is given by
\begin{equation}
\phi_0(x)=C_1e^{\sqrt{2\alpha}x}+C_2e^{-\sqrt{2\alpha}x},
\end{equation}
where $C_1$ and $C_2$ are constants.\\

Since $\phi_{\alpha}$ is bounded as $x$ goes to $-\infty$ so we must have
$C_2=0$.

Using the continuity of $\phi_{\alpha}$ at $x=x_0$ we have $\phi_0(x_0)=x_0$ then
\begin{equation}
C_1e^{\sqrt{2\alpha}x_0}=x_0.
\end{equation}

 Hence $C_1=x_0 e^{-\sqrt{2\alpha}x_0}$.

Then
\begin{equation}
\phi_{\alpha}(t,x)=e^{-\alpha t}x_0 e^{-\sqrt{2\alpha}x_0}e^{\sqrt{2\alpha}x}                                  .
\end{equation}

Using now the high contact equation i.e,  $\phi_{\alpha}$ is $C^1$ at $x=x_0$ we get the following equation
\begin{equation}
x_0 e^{-\sqrt{2\alpha}x_0}e^{\sqrt{2\alpha}x_0}(\sqrt{2\alpha})=1.
\end{equation}
Then we deduce that
\begin{equation}
x_0 =\frac{1}{\sqrt{2\alpha}}.
\end{equation}

In $D^c=\{x\geq x_0\}, x_0\geq0$ we have  $L\phi_{\alpha}+f=L\phi_{\alpha}=-\alpha e^{-\alpha t} x\leq 0.$\\

(ii)By condition (iii) of Theorem 2.2 we have $\psi(x)=f(x)$ on $D$. i:e $\psi(x)=0$ on $D$.
If $x\geq x_0$ we have $(\alpha-A)\phi(x)=\alpha x$.
Then $\psi$ gets the following expression
\begin{equation}
\psi(x)=\begin{cases}
0 , \quad x<x_0,\\
\alpha x, \quad x\geq x_0.
\end{cases}
\end{equation}
Then
\begin{align}
R_{\alpha}\psi(x)&= \int_{x_0}^{\infty}\psi(y)R_{\alpha}(x,dy)\nonumber\\
&=\int_{x_0}^{\infty}\psi(y)\int_0^{\infty} e^{-\alpha t}P_x(B_t\in dy)dt\nonumber\\
&=\alpha\int_0^{\infty} e^{-\alpha t}(\int_{x_0}^{\infty}y P_x(B_t\in dy)dt.
\end{align}

Using the same computation as in \eqref{ex1R} we get that
\begin{equation}
R_{\alpha}\psi(x)=\alpha\int_0^{+\infty}e^{-\alpha t}t\frac{e^{-\frac{(x_0-x)^2}{2t}}}{\sqrt{2\pi t}}dt+\alpha x \int_0^{+\infty}e^{-\alpha t}\int_{x_0-x}^{+\infty}\frac{e^{-\frac{z^2}{2t}}}{\sqrt{2\pi t}}dz dt.
\end{equation}
We distinguish two cases:
\begin{itemize}
\item 1- $x<x_0$, in this case $x_0-x>0$.\\

Using the same calculus as in the first example we get that 
\begin{equation}\label{eq3.59}
R_{\alpha}\psi(x)=\alpha\int_0^{+\infty}e^{-\alpha t}t\frac{e^{-\frac{(x_0-x)^2}{2t}}}{\sqrt{2\pi t}}dt+\frac{\alpha x}{2}\int_{0}^{+\infty}e^{-\alpha t}(1-\int_0^{\frac{x_0-x}{\sqrt{2t}}}\frac{2e^{-z^2}}{\sqrt{\pi}}dz)dt.
\end{equation}

Or we have 
\begin{equation}\label{eq3.60}
\int_{0}^{+\infty}e^{-\alpha t}(1-\int_0^{\frac{x_0-x}{\sqrt{2t}}}\frac{2e^{-z^2}}{\sqrt{\pi}}dz)dt=\frac{1}{\alpha}e^{-\sqrt{2\alpha}(x_0-x)}
\end{equation}
we derive equation \eqref{eq3.60} first with respect to $x$ then $\alpha$ then we get 
\begin{equation}\label{eq3.61}
\int_0^{+\infty}e^{-\alpha t}t\frac{e^{-\frac{(x_0-x)^2}{2t}}}{\sqrt{2\pi t}}dt=e^{-\sqrt{2\alpha}(x_0-x)}(\frac{1}{2\alpha\sqrt{2\alpha}}+\frac{x_0-x}{2\alpha}).
\end{equation}
Replacing \eqref{eq3.60} and \eqref{eq3.61} in equation \eqref{eq3.59} we get 
\begin{equation}\label{eq3.62}
R_{\alpha}\psi(x)=e^{-\sqrt{2\alpha}(x_0-x)}(\frac{1}{2\sqrt{2\alpha}}+\frac{x_0}{2}).
\end{equation}
\item 2- By Theorem 2.2 we have that outside D i.e 
\begin{equation}\label{eq3.63}
x\geq x_0 ,\quad R_{\alpha}\psi(x)=x.
\end{equation}

Then by the continuity of $R_{\alpha}\psi$ at $x=x_0$ and combining  \eqref{eq3.63} and \eqref{eq3.62} at $x=x_0$, we get 
\begin{equation}
x_0=\frac{1}{\sqrt{2\alpha}}.
\end{equation}
\end{itemize}

\begin{example}
Find $\Phi_{\alpha}(x)$ and $\tau^{*}$ such that
\begin{equation}
\Phi_{\alpha}(x)= \sup_{\tau} E_{x}[e^{-\alpha \tau} B(\tau)] = E_{x}[e^{-\alpha \tau^{*}} B(\tau^*)].
\end{equation}
where $B(t)$ is Brownian motion reflected upwards when $B(t)=0$.
(i)The function $\phi_{\alpha}$ should verify the following PDE
\begin{equation}
\begin{cases}
\frac{\partial \phi_{\alpha}}{\partial t}(t,x)+\frac{1}{2}\frac{\partial^2 \phi_{\alpha}}{\partial x^2}(t,x)=0 \text{ on } D,\\
\phi_{\alpha}(t,x)=e^{-\alpha t} x,\quad x\notin D,\\
\frac{\partial \phi_{\alpha}}{\partial x}(t,0)=0.
\end{cases}
\end{equation}
Put $\phi_{\alpha}(t,x)=\phi_0(x) e^{-\alpha t}$
then $\phi_0$ verifies the following second order differential equation
\begin{equation}\label{secondreflectedbis}
\begin{cases}
\frac{1}{2}\frac{\partial^2 \phi_{0}}{\partial x^2}(x)-\alpha \phi_{0}(x)=0 \text{ on } D\\
\phi_{0}(x)=x,\quad x\notin D.\\
\phi_{0}'(0)=0
\end{cases}
\end{equation}
The general solution of the equation \eqref{secondreflectedbis} is given by
\begin{equation}
\phi_0(x)=C_1e^{\sqrt{2\alpha}x}+C_2e^{-\sqrt{2\alpha}x},
\end{equation}
where $C_1$ and $C_2$ are constants.\\
Since we have $\phi_{0}'(0)=0$ then $C_1=C_2$ and
\begin{equation}
\phi_0(x)=C_1(e^{\sqrt{2\alpha}x}+e^{-\sqrt{2\alpha}x}).
\end{equation}
Using the continuity of $\phi_{\alpha}$ at $x=x_0$ we have $\phi_0(x_0)=x_0$ then
\begin{equation}
C_1=\frac{x_0}{e^{\sqrt{2\alpha}x_0}+e^{-\sqrt{2\alpha}x_0}}.
\end{equation}
Then the solution of \eqref{secondreflectedbis} is given by
\begin{equation}\label{eqrefappbis}
\phi_0(x)=\frac{x_0}{e^{\sqrt{2\alpha}x_0}+e^{-\sqrt{2\alpha}x_0}}(e^{\sqrt{2\alpha}x}+e^{-\sqrt{2\alpha}x}).
\end{equation}
Using the high contact condition at $x=x_0$ we get that
\begin{equation}
\frac{\sqrt{2\alpha}x_0}{e^{\sqrt{2\alpha}x_0}+e^{-\sqrt{2\alpha}x_0}}(e^{\sqrt{2\alpha}x_0}-e^{-\sqrt{2\alpha}x_0})=1.
\end{equation}
The $x_0$ is a solution of the following equation
\begin{equation}
\tanh(\sqrt{2\alpha}x)=\frac{1}{\sqrt{2\alpha}x}.
\end{equation}
This equation has two solutions one negative and one positive.\\
For $x_0>0$ we have  $D=\{ 0\leq x<x_0\}$.
\end{example}

(ii)By Theorem 2.2 we have that outside of $D$ i.e
\begin{equation}\label{express1bis}
x\geq x_0,\quad  R_{\alpha}\psi(x)=x.
\end{equation}

 We have $\psi(x)=f(x)=0$ for $x<x_0$ and for  $x\geq x_0, \psi(x)=(\alpha-A)\phi(x)=\alpha x$.
Then we get that for all $x\in \mathbb{R}_+$

\begin{align}
R_{\alpha}\psi(x)&=\int_0^{\infty}e^{-\alpha t}\int_{x_0}^{\infty}\psi(y)P_x(|B_t|\in dy)dt\\
&=\alpha \int_0^{\infty} e^{-\alpha t}\int_{x_0}^{\infty}y\frac{e^{-\frac{(y-x)^2}{2t}}+e^{-\frac{(y+x)^2}{2t}}}{\sqrt{2\pi t}}dy dt\\\label{eqfireqbis}
&=\alpha \int_0^{\infty} \frac{e^{-\alpha t}}{\sqrt{2\pi t}}\int_{x_0}^{\infty}\{ye^{-\frac{(y-x)^2}{2t}}+ye^{-\frac{(y+x)^2}{2t}}\}dy dt\\
&=\alpha \int_0^{\infty} \frac{e^{-\alpha t}}{\sqrt{2\pi
t}}\{-t[e^{-\frac{(y-x)^2}{2t}}]_{x_0}^{\infty}+x\int_{x_0}^{\infty}e^{-\frac{(y-x)^2}{2t}}dy-t[e^{-\frac{(y+x)^2}{2t}}]_{x_0}^{\infty}-x\int_{x_0}^{\infty}e^{-\frac{(y+x)^2}{2t}}dy\}dt\\
&=\alpha \int_0^{\infty} \frac{e^{-\alpha t}}{\sqrt{2\pi
t}}\{te^{-\frac{(x_0-x)^2}{2t}}+x\int_{x_0}^{\infty}e^{-\frac{(y-x)^2}{2t}}dy+te^{-\frac{(x_0+x)^2}{2t}}-x\int_{x_0}^{\infty}e^{-\frac{(y+x)^2}{2t}}dy\}dt\\
&=\alpha \int_0^{\infty} \frac{e^{-\alpha t}}{\sqrt{2\pi
t}}\{te^{-\frac{(x_0-x)^2}{2t}}+te^{-\frac{(x_0+x)^2}{2t}}+x\int_{x_0}^{\infty}e^{-\frac{(y-x)^2}{2t}}dy-x\int_{x_0}^{\infty}e^{-\frac{(y+x)^2}{2t}}dy\}dt.\label{express2bis}
\end{align}

We study now the sum of the two integrals with respect to $dy$ in the previous equation.\\
We have
\begin{align}
&\alpha x\int_0^{+\infty}e^{-\alpha t}\int_{x_0}^{+\infty}\frac{e^{-\frac{(y-x)^2}{2t}}}{\sqrt{2\pi t}}dy dt-\alpha x\int_0^{+\infty}e^{-\alpha t}\int_{x_0}^{+\infty}\frac{e^{-\frac{(y+x)^2}{2t}}}{\sqrt{2\pi t}}dy dt\nonumber\\
&=\alpha x\int_0^{+\infty}e^{-\alpha t}\int_{\frac{x_0-x}{\sqrt{2t}}}^{+\infty}\frac{e^{-z^2}}{\sqrt{\pi }}dz dt-\alpha x\int_0^{+\infty}e^{-\alpha t}\int_{\frac{x_0+x}{\sqrt{2t}}}^{+\infty}\frac{e^{-z^2}}{\sqrt{\pi }}dz dt\nonumber\\
&=\alpha\frac{x}{2}\int_0^{+\infty}e^{-\alpha t}\Big(\int_{\frac{x_0-x}{\sqrt{2t}}}^{+\infty}2\frac{e^{-z^2}}{\sqrt{\pi }}dz-\int_{\frac{x_0+x}{\sqrt{2t}}}^{+\infty}2\frac{e^{-z^2}}{\sqrt{\pi} }dz\Big)dt
\end{align}
For $0\leq  x<x_0$, we have that the last equation is equal to
\begin{align}
&\alpha\frac{x}{2}\int_0^{+\infty}e^{-\alpha t}( 1-H(\frac{\sqrt{2}(x_0-x)}{2\sqrt{t}}))dt-\alpha\frac{x}{2}\int_0^{+\infty}e^{-\alpha t}( 1-H(\frac{\sqrt{2}(x_0+x)}{2\sqrt{t}}))dt\nonumber\\
&=\frac{x}{2}e^{-\sqrt{2\alpha}(x_0-x)}-\frac{x}{2}e^{-\sqrt{2\alpha}(x_0+x)},
\end{align}
where
\begin{equation}
H(u)=\int_0^u\frac{2e^{-z^2}}{\sqrt{\pi}}dz, \quad \forall u>0.
\end{equation}
Then we deduce that for $0\leq x<x_0$
\begin{equation}\label{eqcalc1bis}
R_{\alpha}\psi(x)=\alpha \int_0^{+\infty}\frac{e^{-\alpha t}}{\sqrt{2\pi t}}t(e^{-\frac{(x_0-x)^2}{2t}}+e^{-\frac{(x_0+x)^2}{2t}})dt+\frac{x}{2}(e^{-\sqrt{2\alpha}(x_0-x)}-e^{-\sqrt{2\alpha}(x_0+x)}).
\end{equation}
Or we have
\begin{equation}
\int_0^{+\infty}e^{-\alpha t}\Big(\int_{\frac{x_0-x}{\sqrt{2t}}}^{+\infty}\frac{e^{-z^2}}{\sqrt{\pi }}dz-\int_{\frac{x_0+x}{\sqrt{2t}}}^{+\infty}\frac{e^{-z^2}}{\sqrt{\pi} }dz\Big)dt=\frac{1}{2\alpha}(e^{-\sqrt{2\alpha}(x_0-x)}-e^{-\sqrt{2\alpha}(x_0+x)}).
\end{equation}
We derive the previous equation with respect to $x$, we get
\begin{equation}
\int_0^{+\infty}\frac{e^{-\alpha t}}{\sqrt{2\pi t}}\Big(e^{-\frac{(x_0-x)^2}{2t}}+e^{-\frac{(x_0+x)^2}{2t}}\Big)dt=\frac{1}{\sqrt{2\alpha}}(e^{-\sqrt{2\alpha}(x_0-x)}+e^{-\sqrt{2\alpha}(x_0+x)}).
\end{equation}
We derive now with respect to $\alpha$ we get
\begin{align}\label{eqcalcbis}
&\int_0^{+\infty}t\frac{e^{-\alpha t}}{\sqrt{2\pi t}}\Big(e^{-\frac{(x_0-x)^2}{2t}}+e^{-\frac{(x_0+x)^2}{2t}}\Big)dt=\frac{1}{2\alpha\sqrt{2\alpha}}(e^{-\sqrt{2\alpha}(x_0-x)}+e^{-\sqrt{2\alpha}(x_0+x)})\nonumber\\
&+
\frac{1}{2\alpha}(x_0-x)e^{-\sqrt{2\alpha}(x_0-x)}+\frac{1}{2\alpha}(x_0+x)e^{-\sqrt{2\alpha}(x_0+x)}.
\end{align}

Substituting \eqref{eqcalcbis} in \eqref{eqcalc1bis} we get that
\begin{equation}\label{eqcalc3bis}
R_{\alpha}\psi(x)=(\frac{1}{2\sqrt{2\alpha}}+\frac{1}{2}x_0)e^{-\sqrt{2\alpha}x_0}(e^{\sqrt{2\alpha}x}+e^{-\sqrt{2\alpha}x}), \quad 0\leq x<x_0.
\end{equation}
Let us now study $R_{\alpha}\psi(x_0)$.\\
Using the continuity of $R_{\alpha}\psi$ at $x_0$ i.e combining the two expressions \eqref{express1bis} and \eqref{eqcalc3bis} of $R_{\alpha}\psi$ at $x=x_0$ we get
\begin{equation}\label{eqkhbis}
(\frac{1}{2\sqrt{2\alpha}}+\frac{1}{2}x_0)e^{-\sqrt{2\alpha}x_0}(e^{\sqrt{2\alpha}x_0}+e^{-\sqrt{2\alpha}x_0})=x_0.
\end{equation}
Form this equation we get that
\begin{equation}
(\frac{1}{2\sqrt{2\alpha}}+\frac{1}{2}x_0)e^{-\sqrt{2\alpha}x_0}=\frac{x_0}{e^{\sqrt{2\alpha}x_0}+e^{-\sqrt{2\alpha}x_0}}.
\end{equation}
Replacing this in \eqref{eqcalc3bis} we get
\begin{equation}\label{eqcalc3biss}
R_{\alpha}\psi(x)=\frac{x_0(e^{\sqrt{2\alpha}x}+e^{-\sqrt{2\alpha}x})}{e^{\sqrt{2\alpha}x_0}+e^{-\sqrt{2\alpha}x_0}}, \quad 0\leq x<x_0.
\end{equation}
Developing equation \eqref{eqkhbis} then multiplying it  by $e^{\sqrt{2\alpha}x_0}$, we get

\begin{equation}\label{eqx01bis}
\frac{1}{\sqrt{2\alpha}x_0}(e^{\sqrt{2\alpha}x_0}+e^{-\sqrt{2\alpha}x_0})=e^{\sqrt{2\alpha}x_0}-e^{-\sqrt{2\alpha}x_0}.
\end{equation}
Then we deduce that $x_0$ satisfies the following equation
\begin{equation}\label{eqx0bis}
th(\sqrt{2\alpha}x)=\frac{1}{\sqrt{2\alpha}x}.
\end{equation}
To summarize, we have proved that
\begin{equation}
R_{\alpha}\psi(x)=
\begin{cases}
\frac{x_0(e^{\sqrt{2\alpha}x}+e^{-\sqrt{2\alpha}x})}{e^{\sqrt{2\alpha}x_0}+e^{-\sqrt{2\alpha}x_0}}, \quad 0\leq x<x_0\\
x, \quad x\geq x_0,
\end{cases}
\end{equation}
where $x_0$ satisfies the equation \eqref{eqx0bis}.
The next step is to verify the assertion (ii) of Theorem 2.2 i.e that for all $x\in \mathbb{R}_+$ we have $R_{\alpha}\psi(x)\geq x$.
In fact, consider the function
\begin{equation}
f(x)=R_{\alpha}\psi(x)-x=\begin{cases}
\frac{x_0(e^{\sqrt{2\alpha}x}+e^{-\sqrt{2\alpha}x})}{e^{\sqrt{2\alpha}x_0}+e^{-\sqrt{2\alpha}x_0}}-x, \quad 0\leq x<x_0\\
0, \quad x\geq x_0,
\end{cases}
\end{equation}
For $0\leq x<x_0$, we have
\begin{equation}\label{eqf'bis}
f'(x)=\frac{\sqrt{2\alpha}x_0(e^{\sqrt{2\alpha}x}-e^{-\sqrt{2\alpha}x})}{e^{\sqrt{2\alpha}x_0}+e^{-\sqrt{2\alpha}x_0}}-1.
\end{equation}
Then
\begin{equation}
f'(x)=0\Rightarrow e^{\sqrt{2\alpha}x}-e^{-\sqrt{2\alpha}x}=\frac{e^{\sqrt{2\alpha}x_0}+e^{-\sqrt{2\alpha}x_0}}{\sqrt{2\alpha}x_0}
\end{equation}
From equation \eqref{eqx01bis}, we have that $x_0$ satisfies
\begin{equation}
e^{\sqrt{2\alpha}x_0}-e^{-\sqrt{2\alpha}x_0}=\frac{e^{\sqrt{2\alpha}x_0}+e^{-\sqrt{2\alpha}x_0}}{\sqrt{2\alpha}x_0}.
\end{equation}
Then $x_0$ is the unique positive solution of $f'(x)=0$.\\
Or we have that $f$ is strictly decreasing on $[0,x_0[$. \\
In fact one can write $f'(x)$ in \eqref{eqf'bis} as the following
\begin{equation}
f'(x)=\frac{\sqrt{2\alpha}x_0 sh(\sqrt{2\alpha}x)}{ch(\sqrt{2\alpha}x_0)}-1,
\end{equation}
then
\begin{equation}
f''(x)=\frac{2\alpha ch(\sqrt{2\alpha}x)}{ch(\sqrt{2\alpha}x_0)}>0.
\end{equation}
Then $f'$ is strictly increasing in $[0,x_0[$ with supremum $f'(x_0)=0$. Then we deduce that for $x\in[0,x_0[, f'(x)<0$. Then $f(x)=R_{\alpha}\psi(x)-x$ is decreasing in $[0,x_0[$ with infimum $f(x_0)=R_{\alpha}\psi(x_0)-x_0=0$.
Therefore $f(x)>0$ for all $x\in [0,x_0[$ i.e $R_{\alpha}\psi(x)>x$ for all $x\in [0,x_0[$.\\
In addition for all $x\geq x_0$ we have $R_{\alpha}\psi(x)=x$. Then we deduce that for all $x\in[0,+\infty[, R_{\alpha}\psi(x)\geq x$.
Then the assertions of the verification Theorem 2.2 are well satisfied.

\begin{example}
Find $\Phi_{\alpha}(x)$ and $\tau^{*}$ such that
\begin{equation}
\Phi_{\alpha}(x)= \sup_{\tau} E_{x}[e^{-\alpha \tau} B(\tau)] = E_{x}[e^{-\alpha \tau^{*}} B(\tau^*)].
\end{equation}
where $B(t)$ is Brownian motion trapped when B(t)=0.
More precisely, we consider
\begin{equation}
\tau(0)=\inf\{ t\geq 0,|B(t)=0\}.
\end{equation}
The Brownian motion process trapped at 0 is defined by
\begin{equation}
B_0(t)=B(t\wedge \tau_0).
\end{equation}
From now on we denote $B_0(t)$ by $B(t)$.\\
(i)The function $\phi_{\alpha}$ should verify the following PDE
\begin{equation}
\begin{cases}
\frac{\partial \phi_{\alpha}}{\partial t}(t,x)+\frac{1}{2}\frac{\partial^2 \phi_{\alpha}}{\partial x^2}(t,x)=0 \text{ on } D\\
\phi_{\alpha}(t,x)=e^{-\alpha t} x,\quad x\notin D.\\
\frac{\partial^2 \phi_{\alpha}}{\partial x^2}(t,0)=0
\end{cases}
\end{equation}
Put $\phi_{\alpha}(t,x)=\phi_0(x) e^{-\alpha t}$
then $\phi_0$ verifies the following second order differential equation
\begin{equation}\label{secondreflected}
\begin{cases}
\frac{1}{2}\frac{\partial^2 \phi_{0}}{\partial x^2}(x)-\alpha \phi_{0}(x)=0 \text{ on } D\\
\phi_{0}(x)=x,\quad x\notin D.\\
\phi_{0}''(0)=0
\end{cases}
\end{equation}
The general solution of the equation \eqref{secondreflected} is given by
\begin{equation}
\phi_0(x)=C_1e^{\sqrt{2\alpha}x}+C_2e^{-\sqrt{2\alpha}x}
\end{equation}
where $C_1$ and $C_2$ are constants.\\
Since we have $\phi_{0}''(0)=0$ then $C_2=-C_1$ and
\begin{equation}
\phi_0(x)=C_1(e^{\sqrt{2\alpha}x}-e^{-\sqrt{2\alpha}x})
\end{equation}
Using the continuity of $\phi_{0}$ at $x=x_0$ we have $\phi_0(x_0)=x_0$ then
\begin{equation}
C_1=\frac{x_0}{e^{\sqrt{2\alpha}x_0}-e^{-\sqrt{2\alpha}x_0}}.
\end{equation}
Then the solution of \eqref{secondreflected} is given by
\begin{equation}
\phi_0(x)=\frac{x_0}{e^{\sqrt{2\alpha}x_0}-e^{-\sqrt{2\alpha}x_0}}(e^{\sqrt{2\alpha}x}-e^{-\sqrt{2\alpha}x}).
\end{equation}
Using the high contact condition at $x=x_0$ we get that
\begin{equation}
\frac{\sqrt{2\alpha}x_0}{e^{\sqrt{2\alpha}x_0}-e^{-\sqrt{2\alpha}x_0}}(e^{\sqrt{2\alpha}x_0}+e^{-\sqrt{2\alpha}x_0})=1.
\end{equation}
The $x_0$ is a solution of the following equation
\begin{equation}
\coth(\sqrt{2\alpha}x)=\frac{1}{\sqrt{2\alpha}x}.
\end{equation}
This equation has two solutions one negative and one positive.\\
For $x_0>0$ we have  $D=\{ 0<x<x_0\}$.
\end{example}

(ii)By Theorem 2.2 we have that outside of $D$ i.e
\begin{equation}\label{express1}
x\geq x_0,\quad  R_{\alpha}\psi(x)=x.
\end{equation}

 We have $\psi(x)=f(x)=0$ for $0<x<x_0$ and for  $x\geq x_0, \psi(x)=(\alpha-A)\phi(x)=\alpha x$.\\
 In addition the density of Brownian motion trapped at 0 is given by
 \begin{equation}
 P_x(B_t\in dy)=\frac{e^{-\frac{(y-x)^2}{2t}}-e^{-\frac{(y+x)^2}{2t}}}{\sqrt{2\pi t}}dy.
 \end{equation}
Then we get that for all $x\in \mathbb{R}_+^*$

\begin{align}
R_{\alpha}\psi(x)&=\int_0^{\infty}e^{-\alpha t}\int_{x_0}^{\infty}\psi(y)P_x(B_t\in dy)dt\\
&=\alpha \int_0^{\infty} e^{-\alpha t}\int_{x_0}^{\infty}y\frac{e^{-\frac{(y-x)^2}{2t}}-e^{-\frac{(y+x)^2}{2t}}}{\sqrt{2\pi t}}dy dt\\
&=\alpha \int_0^{\infty} \frac{e^{-\alpha t}}{\sqrt{2\pi t}}\int_{x_0}^{\infty}\{ye^{-\frac{(y-x)^2}{2t}}-ye^{-\frac{(y+x)^2}{2t}}\}dy dt\\
&=\alpha \int_0^{\infty} \frac{e^{-\alpha t}}{\sqrt{2\pi
t}}\{-t[e^{-\frac{(y-x)^2}{2t}}]_{x_0}^{\infty}+x\int_{x_0}^{\infty}e^{-\frac{(y-x)^2}{2t}}dy+t[e^{-\frac{(y+x)^2}{2t}}]_{x_0}^{\infty}+x\int_{x_0}^{\infty}e^{-\frac{(y+x)^2}{2t}}dy\}dt\\
&=\alpha \int_0^{\infty} \frac{e^{-\alpha t}}{\sqrt{2\pi
t}}\{te^{-\frac{(x_0-x)^2}{2t}}+x\int_{x_0}^{\infty}e^{-\frac{(y-x)^2}{2t}}dy-te^{-\frac{(x_0+x)^2}{2t}}+x\int_{x_0}^{\infty}e^{-\frac{(y+x)^2}{2t}}dy\}dt\\
&=\alpha \int_0^{\infty} \frac{e^{-\alpha t}}{\sqrt{2\pi
t}}\{te^{-\frac{(x_0-x)^2}{2t}}-te^{-\frac{(x_0+x)^2}{2t}}+x\int_{x_0}^{\infty}e^{-\frac{(y-x)^2}{2t}}dy+x\int_{x_0}^{\infty}e^{-\frac{(y+x)^2}{2t}}dy\}dt.\label{express2}
\end{align}

We study now the sum of the two integrals with respect to $dy$ in the previous equation.\\
We have
\begin{align}
&\alpha x\int_0^{+\infty}e^{-\alpha t}\int_{x_0}^{+\infty}\frac{e^{-\frac{(y-x)^2}{2t}}}{\sqrt{2\pi t}}dy dt+\alpha x\int_0^{+\infty}e^{-\alpha t}\int_{x_0}^{+\infty}\frac{e^{-\frac{(y+x)^2}{2t}}}{\sqrt{2\pi t}}dy dt\nonumber\\
&=\alpha x\int_0^{+\infty}e^{-\alpha t}\int_{\frac{x_0-x}{\sqrt{2t}}}^{+\infty}\frac{e^{-z^2}}{\sqrt{\pi }}dz dt+\alpha x\int_0^{+\infty}e^{-\alpha t}\int_{\frac{x_0+x}{\sqrt{2t}}}^{+\infty}\frac{e^{-z^2}}{\sqrt{\pi }}dz dt\nonumber\\
&=\alpha\frac{x}{2}\int_0^{+\infty}e^{-\alpha t}\Big(\int_{\frac{x_0-x}{\sqrt{2t}}}^{+\infty}2\frac{e^{-z^2}}{\sqrt{\pi }}dz+\int_{\frac{x_0+x}{\sqrt{2t}}}^{+\infty}2\frac{e^{-z^2}}{\sqrt{\pi} }dz\Big)dt
\end{align}
For $0<x<x_0$, we have that the last equation is equal to
\begin{align}
&=\alpha\frac{x}{2}\int_0^{+\infty}e^{-\alpha t}( 1-H(\frac{\sqrt{2}(x_0-x)}{2\sqrt{t}}))dt+\alpha\frac{x}{2}\int_0^{+\infty}e^{-\alpha t}( 1-H(\frac{\sqrt{2}(x_0+x)}{2\sqrt{t}}))dt\nonumber\\
&=\frac{x}{2}e^{-\sqrt{2\alpha}(x_0-x)}+\frac{x}{2}e^{-\sqrt{2\alpha}(x_0+x)},
\end{align}
where
\begin{equation}
H(u)=\int_0^u\frac{2e^{-z^2}}{\sqrt{\pi}}dz, \quad \forall u>0.
\end{equation}
Then we deduce that for $0< x<x_0$
\begin{equation}\label{eqcalc1}
R_{\alpha}\psi(x)=\alpha \int_0^{+\infty}\frac{e^{-\alpha t}}{\sqrt{2\pi t}}t(e^{-\frac{(x_0-x)^2}{2t}}-e^{-\frac{(x_0+x)^2}{2t}})dt+\frac{x}{2}(e^{-\sqrt{2\alpha}(x_0-x)}+e^{-\sqrt{2\alpha}(x_0+x)}).
\end{equation}
Or we have
\begin{equation}
\int_0^{+\infty}e^{-\alpha t}\Big(\int_{\frac{x_0-x}{\sqrt{2t}}}^{+\infty}\frac{e^{-z^2}}{\sqrt{\pi }}dz+\int_{\frac{x_0+x}{\sqrt{2t}}}^{+\infty}\frac{e^{-z^2}}{\sqrt{\pi} }dz\Big)dt=\frac{1}{2\alpha}(e^{-\sqrt{2\alpha}(x_0-x)}+e^{-\sqrt{2\alpha}(x_0+x)}).
\end{equation}
We derive the previous equation with respect to $x$, we get
\begin{equation}
\int_0^{+\infty}\frac{e^{-\alpha t}}{\sqrt{2\pi t}}\Big(e^{-\frac{(x_0-x)^2}{2t}}+e^{-\frac{(x_0+x)^2}{2t}}\Big)dt=\frac{1}{\sqrt{2\alpha}}(e^{-\sqrt{2\alpha}(x_0-x)}+e^{-\sqrt{2\alpha}(x_0+x)}).
\end{equation}
We derive now with respect to $\alpha$ we get
\begin{align}\label{eqcalc}
&\int_0^{+\infty}t\frac{e^{-\alpha t}}{\sqrt{2\pi t}}\Big(e^{-\frac{(x_0-x)^2}{2t}}-e^{-\frac{(x_0+x)^2}{2t}}\Big)dt=\frac{1}{2\alpha\sqrt{2\alpha}}(e^{-\sqrt{2\alpha}(x_0-x)}-e^{-\sqrt{2\alpha}(x_0+x)})\nonumber\\
&+
\frac{1}{2\alpha}(x_0-x)e^{-\sqrt{2\alpha}(x_0-x)}-\frac{1}{2\alpha}(x_0+x)e^{-\sqrt{2\alpha}(x_0+x)}.
\end{align}

Substituting \eqref{eqcalc} in \eqref{eqcalc1} we get that
\begin{equation}\label{eqcalc3}
R_{\alpha}\psi(x)=(\frac{1}{2\sqrt{2\alpha}}+\frac{1}{2}x_0)e^{-\sqrt{2\alpha}x_0}(e^{\sqrt{2\alpha}x}-e^{-\sqrt{2\alpha}x}), \quad 0<x<x_0.
\end{equation}
Let us now study $R_{\alpha}\psi(x_0)$.\\
Using the continuity of $R_{\alpha}\psi$ at $x_0$ i.e combining the two expressions \eqref{express1} and \eqref{eqcalc3} of $R_{\alpha}\psi$ at $x=x_0$ we get
\begin{equation}\label{eqkh}
(\frac{1}{2\sqrt{2\alpha}}+\frac{1}{2}x_0)e^{-\sqrt{2\alpha}x_0}(e^{\sqrt{2\alpha}x_0}-e^{-\sqrt{2\alpha}x_0})=x_0.
\end{equation}
Form this equation we get that
\begin{equation}
(\frac{1}{2\sqrt{2\alpha}}+\frac{1}{2}x_0)e^{-\sqrt{2\alpha}x_0}=\frac{x_0}{e^{\sqrt{2\alpha}x_0}-e^{-\sqrt{2\alpha}x_0}}.
\end{equation}
Replacing this in \eqref{eqcalc3} we get
\begin{equation}\label{eqcalc3bisss }
R_{\alpha}\psi(x)=\frac{x_0(e^{\sqrt{2\alpha}x}-e^{-\sqrt{2\alpha}x})}{e^{\sqrt{2\alpha}x_0}-e^{-\sqrt{2\alpha}x_0}}, \quad 0<x<x_0.
\end{equation}
Developing equation \eqref{eqkh} then multiplying it  by $e^{\sqrt{2\alpha}x_0}$, we get

\begin{equation}\label{eqx01}
\frac{1}{\sqrt{2\alpha}x_0}(e^{\sqrt{2\alpha}x_0}-e^{-\sqrt{2\alpha}x_0})=e^{\sqrt{2\alpha}x_0}+e^{-\sqrt{2\alpha}x_0}.
\end{equation}
Then we deduce that $x_0$ satisfies the following equation
\begin{equation}\label{eqx0}
coth(\sqrt{2\alpha}x)=\frac{1}{\sqrt{2\alpha}x}.
\end{equation}
To summarize, we have proved that
\begin{equation}
R_{\alpha}\psi(x)=
\begin{cases}
\frac{x_0(e^{\sqrt{2\alpha}x}-e^{-\sqrt{2\alpha}x})}{e^{\sqrt{2\alpha}x_0}-e^{-\sqrt{2\alpha}x_0}}, \quad 0<x<x_0\\
x, \quad x\geq x_0,
\end{cases}
\end{equation}
where $x_0$ satisfies the equation \eqref{eqx0}.
The next step is to verify the assertion (ii) of Theorem 2.2 i.e that for all $x\in \mathbb{R}_+$ we have $R_{\alpha}\psi(x)\geq x$.
In fact, consider the function
\begin{equation}
f(x)=R_{\alpha}\psi(x)-x=\begin{cases}
\frac{x_0(e^{\sqrt{2\alpha}x}-e^{-\sqrt{2\alpha}x})}{e^{\sqrt{2\alpha}x_0}-    e^{-\sqrt{2\alpha}x_0}}-x, \quad 0< x<x_0\\
0, \quad x\geq x_0,
\end{cases}
\end{equation}
For $0< x<x_0$, we have
\begin{equation}\label{eqf'}
f'(x)=\frac{\sqrt{2\alpha}x_0(e^{\sqrt{2\alpha}x}+e^{-\sqrt{2\alpha}x})}{e^{\sqrt{2\alpha}x_0}-e^{-\sqrt{2\alpha}x_0}}-1.
\end{equation}
Then
\begin{equation}
f'(x)=0\Rightarrow e^{\sqrt{2\alpha}x}+
e^{-\sqrt{2\alpha}x}=\frac{e^{\sqrt{2\alpha}x_0}-e^{-\sqrt{2\alpha}x_0}}{\sqrt{2\alpha}x_0}
\end{equation}
\begin{figure}
	\centering
	\includegraphics[width=0.9\linewidth]{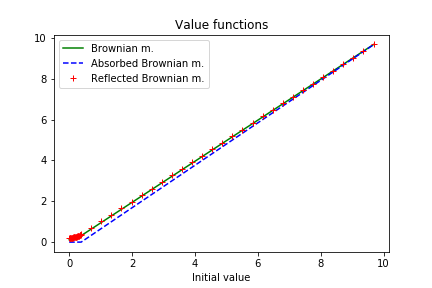}
	\caption{\small{}}
	\label{}
\end{figure}
\begin{figure}
	\centering
	\includegraphics[width=.7\linewidth]{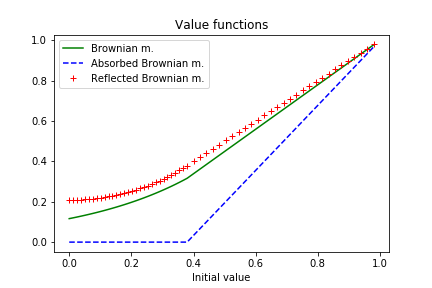}
	\caption{\small{}}
	\label{}
\end{figure}
From equation \eqref{eqx01}, we have that $x_0$ satisfies
\begin{equation}
e^{\sqrt{2\alpha}x_0}+
e^{-\sqrt{2\alpha}x_0}=\frac{e^{\sqrt{2\alpha}x_0}-e^{-\sqrt{2\alpha}x_0}}{\sqrt{2\alpha}x_0}.
\end{equation}
Then $x_0$ is the unique positive solution of $f'(x)=0$.\\
Or we have that $f$ is strictly decreasing on $]0,x_0[$. \\
In fact one can write $f'(x)$ in \eqref{eqf'} as the following
\begin{equation}
f'(x)=\frac{\sqrt{2\alpha}x_0 ch(\sqrt{2\alpha}x)}{sh(\sqrt{2\alpha}x_0)}-1,
\end{equation}
then
\begin{equation}
f''(x)=\frac{2\alpha sh(\sqrt{2\alpha}x)}{sh(\sqrt{2\alpha}x_0)}>0\quad \forall x\in]0,x_0[.
\end{equation}
Then $f'$ is strictly increasing in $]0,x_0[$ with supremum $f'(x_0)=0$. Then we deduce that for $x\in]0,x_0[, f'(x)<0$. Then $f(x)=R_{\alpha}\psi(x)-x$ is strictly decreasing in $]0,x_0[$ with infimum  $f(x_0)=R_{\alpha}\psi(x_0)-x_0=0$.
Therefore $f(x)> 0$ for all $x\in [0,x_0[$ i.e $R_{\alpha}\psi(x)\geq x$ for all $x\in ]0,x_0[$.\\
In addition for all $x\geq x_0$ we have $R_{\alpha}\psi(x)=x$. Then we deduce that for all $x\in]0,+\infty[, R_{\alpha}\psi(x)\geq x$.
Then the assertions of the verification Theorem 2.2 are well satisfied.
\begin{remark}
We denote by $x_0^o,x_0^r$ and $x_0^a$ the values of $x_0$ associated to  the optimal stopping barriers for the usual Brownian motion i.e Example 2.2, reflected Brownian motion at 0 i.e Example 2.3 and absorbed (trapped) Brownian motion i.e Example 2.4 , respectively.
Comparing these three values we have
\begin{equation}
x_0^a<x_0^o<x_0^r.
\end{equation}

Indeed this result is expected because if $B(t)$ is absorbed at 0 then the payoff is 0, and this is the worst that can happen. Therefore one is afraid to wait for large value of $B(t)$ before stopping.
For the reflected Brownian motion, however, there is no disaster if it hits 0 because it is just reflected back. Therefore one can wait for a large value of $B(t)$ before stopping.
\end{remark}

\end{document}